\documentclass[12pt]{article}
\usepackage{amsfonts}
\usepackage{mathrsfs}
\usepackage{amsmath,amssymb}
\baselineskip 5.5pt \lineskip 5.5pt \numberwithin{equation}{section}

\def\span{{\rm span}}

\def\adddot{$\!\!\!${\bf.}\ \ }
\def\adddot{}
\def\ad{\mbox{ad}}
\def\QED{\hfill$\Box$\par}
\def\a{\alpha}

\def\b{\beta}

\def\D{\Delta}

\def\SS{\mathcal {S}}

\def\l{\lambda}

\def\Vir{\mbox{$\mathcal {V}ir$}}

\def\LL{\mathcal {L}}

\def\ssc{\scriptscriptstyle}

\def\cl{\centerline}

\def\lll{\leftline}

\def\D{\Delta}

\def\vs{\vspace*}

\def\H{{\cal H}}

\def\C{\mathbb{C}}

\def\F{\mathbb{F}}
\def\Z{\mathbb{Z}}

\def\N{\mathbb{N}}

\newtheorem{case}{Case}

\newtheorem{theo}{Theorem}[section]

\newtheorem{lemm}[theo]{Lemma}

\newtheorem{rema}[theo]{Remark}
\newtheorem{defi}[theo]{Definition}

\def\ma{\mathbb}
\def\vs{\vspace*}
\def\BE {\begin{eqnarray}}
\def\EE {\end{eqnarray}}
\def\BC {\begin{eqnarray*}}
\def\EC {\end{eqnarray*}}

\def\vs{\vspace*}\def\cl{\centerline}
\begin{document}

\cl{{\bf\large Representations of the Schr\"{o}dinger-Virasoro
algebras}\footnote {Supported by NSF grants 10471091, 10671027 of
China, ``One Hundred Talents Program'' from University of Science
and Technology of China.\\[2pt] \indent Corresponding E-mail:
sd\_junbo@163.com}} \vs{6pt}

\cl{Junbo Li$^{*,\dag)}$, Yucai Su$^{\ddag)}$}

\cl{\small $^{*)}$Department of Mathematics, Shanghai Jiao Tong
University, Shanghai 200240, China}

\cl{\small $^{\dag)}$Department of Mathematics, Changshu Institute
of Technology, Changshu 215500, China}

\cl{\small $^{\ddag)}$Department of Mathematics, University of
Science and Technology of China, Hefei 230026, China}

\cl{\small E-mail: sd\_junbo@163.com, ycsu@ustc.edu.cn} \vs{6pt}

\noindent{\bf{Abstract.}} In this paper it is proved that an
irreducible weight module with finite-dimensional weight spaces over
the  Schr\"{o}dinger-Virasoro algebras is a highest/lowest weight
module or a uniformly bounded module. Furthermore, indecomposable
modules of the intermediate series over these algebras are
completely determined.


\noindent{\it{MR(2000) Subject Classification}: 17B10, 17B65,
17B68.}\vs{12pt}

\lll{\bf1. \ Introduction}
\setcounter{section}{1}\setcounter{theo}{0} \setcounter{equation}{0}
Let $s=0$ or $\frac12$\,. The {\it Schr\"{o}dinger-Virasoro algebra}
$\LL[s]$ introduced in \cite{H1,H2,HU}, in the context of
non-equilibrium statistical physics as a by-product of the
computation of $n$-point functions that are covariant under the
action of the Schr\"{o}dinger group, is the infinite-dimensional Lie
algebra with $\C$-basis $\{L_m,\,Y_p\,,\,M_n,\,c\,|\,m,\,n\in
\Z,\,p\in s+\Z\}$ and Lie brackets,
\begin{eqnarray}
\!\!\!&\!\!\!&
[L_m,\,L_{m'}]=(m'-m)L_{m'+m}+\delta_{m,-m'}\frac{m^3-m}{12}c,\label{1eq1.1}
\\
\!\!\!&\!\!\!&[L_m,\,Y_p]=(p-\frac{m}{2})Y_{p+m},\label{1eq1.2}\\[2pt]
\!\!\!&\!\!\!& [L_m,\,M_n]=nM_{n+m},\label{1eq1.3}\\[3pt]
\!\!\!&\!\!\!&[Y_p,\,Y_{p\,'}]=(p\,'-p)M_{p\,'+p}\,,\label{1eq1.4}\\[3pt]
\!\!\!&\!\!\!&[Y_p,\,M_n]=[M_n,\,M_{n'}]=[\LL,\,c]=0.\label{1eq1.5}
\end{eqnarray}
The Lie algebra $\LL[s]$ contains as subalgebras both the Lie
algebra $\SS$ of invariance of the free Schr\"{o}dinger equation
 and the central
charge-free Virasoro algebra $\Vir$, where $\SS$ is the infinite
dimensional Lie algebra, called the {\it Schr\"{o}dinger algebra},
 with the $\C$-basis
$\{Y_p\,,\,M_n\,|\,n\in\Z,\,p\in s+\Z\}$, and $\Vir$ is the {\it
Virasoro algebra} with the $\C$-basis $\{L_n,\,c\,|\,n\in \Z\}$. The
Lie algebra $\LL[\frac12]$ is called the {\it original
Schr\"{o}din\-ger-Virasoro algebra}, while $\LL[0]$ is called the
{\it twisted Schr\"{o}din\-ger-Virasoro algebra}\,.

 It is well known
that one attempt of introducing two-dimensional conformal field
theory is to understand the universal behavior of two-dimensional
statistical systems at equilibrium and at the critical temperature.
A systematic investigation of the theory of representations of the
Virasoro algebra in the 80's led to introduce the unitary minimal
models, corresponding to the unitary highest weight representations
of the Virasoro algebra with central charge less than one.
Miraculously, covariance alone is enough to allow the computation of
the $n$-point functions for these highly constrained physical
models. Since  both original and twisted Schr\"{o}dinger-Virasoro
Lie algebras are closely related to the Schr\"{o}dinger Lie algebra
and the Virasoro Lie algebra, it is highly expected that they should
consequently play an important role akin to that of the Virasoro Lie
algebra in two-dimensional equilibrium statistical physics (see,
e.g., \cite{H1,H2,HU,LS,RU,U}).

Partially due to the above-stated reasons, the
Schr\"{o}dinger-Virasoro Lie algebras have recently drawn some
attentions in the literature. In particular,  the sets of generators
provided by the cohomology classes of the cocycles for both original
and twisted Schr\"{o}dinger-Virasoro Lie algebras were presented in
\cite{RU}, and the derivation algebra and the automorphism group of
the twisted sector were determined in \cite{LS}. Furthermore, vertex
algebra representations of the Schr\"{o}dinger-Virasoro Lie algebras
out of a charged symplectic boson and a free boson with its
associated vertex operators were constructed in \cite{U}.

Our motivation in studying the Schr\"{o}dinger-Virasoro Lie algebras
is to have a better understanding of their representations. Let us
formulate our main results below.

Denote
\begin{equation}\label{H===}\H=\span\{L_0,M_0,c\},\end{equation} which
is a maximal torus of $\LL[s]$ (note that in case of $s=0$, since
the adjoint operator $\ad_{Y_0}$ is not semi-simple, $\H\oplus\C
Y_0$ is not a torus but a Cartan subalgebra of $\LL[0]$). Denote by
$\H^*$ the dual space of $\H$.
\begin{defi}\adddot\label{defi1-}
\rm A module $V$ over $\LL[s]$ is called a
\begin{enumerate} \item[(i)] {\it
Harish-Chandra module} if $V$ admits a finite-dimensional weight
space decomposition $V=\oplus_{\l\in\H^*} V^\l$, where
$V^\lambda=\{v\in V\,|\,x v=\lambda(x)v,\,x\in{\cal H} \}$ such that
${\rm dim\ssc\,}V^\l<\infty$ for all $\l\in \H^*$ (in case
$V^\l\ne0$, we call $\l$ a {\it weight of $V$});
\item[(ii)] {\it uniformly bounded module} if it is a Harish-Chandra
module such that there exists some $N>0$ with ${\rm
dim\ssc\,}V^\lambda\leqslant N$ for all $\lambda\in{\cal H}^*$;
\item[(iii)] {\it module
of the intermediate series} if $V$ is an indecomposable
Harish-Chandra module  such that ${\rm dim\ssc\,}V^\lambda\leqslant
1$ for all $\lambda\in{\cal H}^*$.
\end{enumerate}\end{defi}
\begin{rema}\rm\begin{enumerate}\item[(i)]
Since the Schr\"{o}dinger subalgebra $\SS$ is an ideal of ${\cal
L}$, if $\SS$ acts trivially on a module $V$, then $V$ is simply a
module over $\mathcal {V}ir$. Thus in the following, we always
suppose that $\SS$ acts nontrivially on $V$.
\item[(ii)]
It is well known (see, e.g., \cite{S3,S4}) that a module $V$ of the
intermediate series over $\Vir$ is a quotient of one of the modules
$A_{a,b},\,A(\a),\,B(\a)$ for some $a,b,\a\in\C$, they all have the
basis $\{x_k\,|\,k\in\Z\}$ such that $c$ acts trivially, and for
$n,k\in\Z$,
\begin{eqnarray}
A_{a,b}:\,\,L_nx_k\!\!\!&=&\!\!\!(a+k+bn)x_{k+n},\label{condition+}\\
A(\a):\,\,L_nx_k\!\!\!&=&\!\!\!(k+n)x_{k+n}\ \,(k\neq0),\ \ \
L_nx_0=n(n+\a)x_n,\label{condition++}\\
B(\a):\,\,L_nx_k\!\!\!&=&\!\!\!kx_{k+n}\ \,(k\neq-n),\ \ \ \ \
 \ \ \ L_nx_{-n}=-n(n+\a)x_0.\label{condition+++}
\end{eqnarray}
\end{enumerate}\end{rema}

The main results of this paper are presented in the following
theorem.
\begin{theo}\adddot\label{maintheo}
\begin{enumerate}\item[\rm(i)]
An irreducible Harish-Chandra module over $\LL[s]$ is either a
highest/lowest weight module or a uniformly bounded one.
\item[\rm(ii)] A module of the intermediate series over $\LL[0]$
is simply a module of the intermediate series over $\Vir$.
\item[\rm(iii)] A module $V$ of the intermediate series
over $\LL[\frac12]$ such that $\SS$ acts nontrivially on $V$ is one
of the modules $$A_{a,b},\ \ B_{a,b},\ \ C_{a},\ \ D_{a},\ \
A_1(\a),\ \ A_2({\a}),\ \ B_1(\a),\ \ B_2({\a}),\ \ C(\a,\a'),\ \
D(\b,\b'),$$ or one of their quotients for
$a,\,b,\,\a,\,\a',\,\beta,\,\b'\in{\ma C}$, whose module structures
are given as follows $($the central element $c$ acts as zero$)$,
where $k\in\frac12\Z,\,i,n\in\Z,\,j,p\in\frac12+\Z$,
\begin{eqnarray}
\!\!\!\!\!\!\!\!\!\!\!\!\!\!\!\!\!\!\!\!A_{a,b}:&&\!\!\!\!\!\!\!M_nx_k=0,\label{caseIIact1}\\
\!\!\!\!\!\!\!\!\!\!\!\!\!\!\!\!\!\!\!\!&&\!\!\!\!\!\!\!Y_px_i=(a+i+2b p)x_{i+p},\ \ Y_px_j=0,\label{caseIIact2}\\
\!\!\!\!\!\!\!\!\!\!\!\!\!\!\!\!\!\!\!\!&&\!\!\!\!\!\!\!L_nx_i=(a+i+bn)x_{i+n},\
\
L_nx_j=\big(a+j+(b+\frac{1}{2})n\big)x_{j+n},\label{caseIIact3}\\
\!\!\!\!\!\!\!\!\!\!\!\!\!\!\!\!\!\!\!\!B_{a,b}:&&\!\!\!\!\!\!\! M_nx_k=0,\label{caseIIIac1}\\
\!\!\!\!\!\!\!\!\!\!\!\!\!\!\!\!\!\!\!\!&&\!\!\!\!\!\!\!Y_px_i=0,\ \ Y_px_j= x_{j+p},\label{caseIIIac2}\\
\!\!\!\!\!\!\!\!\!\!\!\!\!\!\!\!\!\!\!\!&&\!\!\!\!\!\!\!L_nx_i=
(a+i+bn)x_{i+n},\ \ L_nx_j=
\big(a+j+(b+\frac{1}{2})n\big)x_{j+n},\label{caseIIIac3}\\
\!\!\!\!\!\!\!\!\!\!\!\!\!\!\!\!\!\!\!\!C_a:&&\!\!\!\!\!\!\! M_nx_k=0,\label{eqaCa1}\\
\!\!\!\!\!\!\!\!\!\!\!\!\!\!\!\!\!\!\!\!&&\!\!\!\!\!\!\!Y_px_i=(a+i)(a+i+2p)x_{i+p},\ \ Y_px_j=0,\label{eqaCa2}\\
\!\!\!\!\!\!\!\!\!\!\!\!\!\!\!\!\!\!\!\!&&\!\!\!\!\!\!\!L_nx_i=(a+i)x_{i+n},\ \ L_nx_j=(a+j+\frac{3}{2}n)x_{j+n},\label{eqaCa3}\\
\!\!\!\!\!\!\!\!\!\!\!\!\!\!\!\!\!\!\!\!D_a:&&\!\!\!\!\!\!\! M_nx_k=0,\label{eqaDa1}\\
\!\!\!\!\!\!\!\!\!\!\!\!\!\!\!\!\!\!\!\!&&\!\!\!\!\!\!\!Y_px_i= (a+i+p)(a+i-p)x_{i+n},\ \ Y_px_j= 0,\label{eqaDa2}\\
\!\!\!\!\!\!\!\!\!\!\!\!\!\!\!\!\!\!\!\!&&\!\!\!\!\!\!\!L_nx_i=
(a+i-\frac{1}{2}n)x_{i+n},\ \ L_nx_j=(a+j+n)x_{j+n},\label{eqaDa3}
\\ 
\!\!\!\!\!\!\!\!\!\!\!\!\!\!\!\!\!\!\!\!A_1(\a):&&\!\!\!\!\!\!\!
M_nx_k=0,
\label{3caseIIact122}\\
\!\!\!\!\!\!\!\!\!\!\!\!\!\!\!\!\!\!\!\!&&\!\!\!\!\!\!\!Y_px_i=
(-\frac{1}{2}+i+p)x_{i+p},\ \ Y_px_j=0,\ \ L_nx_{\frac12}= n(n+\a)x_{\frac{1}{2}+n},\label{3caseIIact222}\\
\!\!\!\!\!\!\!\!\!\!\!\!\!\!\!\!\!\!\!\!&&\!\!\!\!\!\!\!
L_nx_i=(-\frac{1}{2}+i+\frac{1}{2}n)x_{i+n},\ \ L_nx_j=
(-\frac{1}{2}+j+n)x_{j+n} \ (j\ne\frac12), \label{3caseIIact322}
\end{eqnarray}\begin{eqnarray}
\!\!\!\!\!\!\!\!\!\!\!\!\!\!\!\!\!\!\!\!A_2(\a):&&\!\!\!\!\!\!\!M_nx_k=0,
\label{2caseIIsub2act1}\\
\!\!\!\!\!\!\!\!\!\!\!\!\!\!\!\!\!\!\!\!&&\!\!\!\!\!\!\!
Y_px_i=ix_{i+p},\ \ Y_px_j= 0,\ \ L_nx_{-n}=
-n(n+\a)x_{0},\label{2caseIIsub2act2}
\\
\!\!\!\!\!\!\!\!\!\!\!\!\!\!\!\!\!\!\!\!&&\!\!\!\!\!\!\!
L_nx_i=ix_{i+n}\ (i\neq-n),\ \
L_nx_j=(j+\frac{1}{2}n)x_{j+n},\label{2caseIIsub2act3}
\\
\!\!\!\!\!\!\!\!\!\!\!\!\!\!\!\!\!\!\!\!B_1(\a):&&\!\!\!\!\!\!\! M_nx_k=0,\label{caseIIIac1a}\\
\!\!\!\!\!\!\!\!\!\!\!\!\!\!\!\!\!\!\!\!&&\!\!\!\!\!\!\!Y_px_i=0,\ \
Y_px_j=x_{j+p},\ \ L_nx_0=n(n+\a)x_{n},\label{caseIIIac2a}\\
\!\!\!\!\!\!\!\!\!\!\!\!\!\!\!\!\!\!\!\!&&\!\!\!\!\!\!\!L_nx_i=(i+n)x_{i+n}\ (i\ne0),\ \ L_nx_j=(j+\frac{3}{2}n)x_{j+n},\label{caseIIIac3a}\\
\!\!\!\!\!\!\!\!\!\!\!\!\!\!\!\!\!\!\!\!B_2(\a):&&\!\!\!\!\!\!\! M_nx_k=0,\label{caseIIIac1b}\\
\!\!\!\!\!\!\!\!\!\!\!\!\!\!\!\!\!\!\!\!&&\!\!\!\!\!\!\!Y_px_i=0,\ \
Y_px_j= x_{j+p},\ \ L_nx_{\frac12-n}=-n(n+\a)x_{\frac{1}{2}},\label{caseIIIac2b}\\
\!\!\!\!\!\!\!\!\!\!\!\!\!\!\!\!\!\!\!\!&&\!\!\!\!\!\!\!L_nx_i=
(-\frac{1}{2}+i-\frac{1}{2}n)x_{i+n}, \ \
L_nx_j=(-\frac{1}{2}+j)x_{j+n}\ (j\ne\frac{1}{2}-n),\label{caseIIIac3b}\\
\!\!\!\!\!\!\!\!\!\!\!\!\!\!\!\!\!\!\!\!C(\a,\a'):&&\!\!\!\!\!\!\!M_nx_{-n}
=-2n\a'x_0,\ M_nx_k=0\ (k\ne-n),\label{eqacIIS1+}\\
\!\!\!\!\!\!\!\!\!\!\!\!\!\!\!\!\!\!\!\!&&\!\!\!\!\!\!\!Y_px_i=
i(i+2p)x_{i+p},\ \ Y_px_{-p}=\a'x_0,\ \ Y_px_j=0\ (j\ne-p),\label{eqacIIS2}\\
\!\!\!\!\!\!\!\!\!\!\!\!\!\!\!\!\!\!\!\!&&\!\!\!\!\!\!\!L_nx_{-n}=-n(n\!+\!\a)x_0,\
L_nx_i=ix_{i+n}\, (i\!\ne\!-n),\  L_nx_j=(j\!+\!\frac{3}{2}n)x_{j+n},\label{eqacIIS3}\\
\!\!\!\!\!\!\!\!\!\!\!\!\!\!\!\!\!\!\!\!D(\b,\b'):
&&\!\!\!\!\!\!\!\!\!M_nx_{\frac12}\!=\!2n\b'x_{n+\frac{1}{2}},\
M_nx_k\!=\!0\, (k\!\neq\!\frac{1}{2}),\ Y_px_{\frac12}\!=\!\b'x_{p+\frac{1}{2}}, \ Y_px_j\!=\!0\, (j\!\ne\!\frac12),\label{eqLLSSI1}\\
\!\!\!\!\!\!\!\!\!\!\!\!\!\!\!\!\!\!\!\!&&\!\!\!\!\!\!\!Y_px_i=
(-\frac{1}{2}+i+p)(-\frac{1}{2}+i-p)x_{i+n}, \ \
L_nx_{\frac12}=n(n+\b)x_{n+\frac{1}{2}},
\label{eqLLSSI2}\\
\!\!\!\!\!\!\!\!\!\!\!\!\!\!\!\!\!\!\!\!&&\!\!\!\!\!\!\!L_nx_i=(-\frac{1}{2}+i-\frac{1}{2}n)x_{i+n},\
\ L_nx_j=(-\frac{1}{2}+j+n)x_{j+n}\
(j\ne\frac{1}{2}).\label{eqLLSSI3}
\end{eqnarray}
\end{enumerate}\end{theo}

Throughout the paper, we denote the set of all nonzero integers by
$\Z^*$ and that of all positive integers by $\N$.

\vskip16pt

\lll{\bf2. \ Proof of Theorem \ref{maintheo}(i)}
\setcounter{section}{2}\setcounter{theo}{0} \setcounter{equation}{0}

The proof of Theorem \ref{maintheo} will be divided by several
lemmas. Let $V$ be an indecomposable module over $\LL[s]$. Since
$M_0$ and $c$ are central elements, whose actions on $V$ must be
constants. Thus from (\ref{H===}), we can simply regard the weight
space as the eigenspace of $L_0$, i.e., $V=\oplus_{\l\in\C} V^\l,$
where $V^\lambda=\{v\in V\,|\,L_0 v=\lambda v \}$ for $\l\in\C.$
\begin{lemm}\adddot\label{lemm2.2---}
Fix an $a\in\C$ such that $V^a\ne0$. We have
$V=\oplus_{n\in\frac12\Z}V_n, $ where $V_n=V^{a+n}=\{v\in
V\,|\,L_0v=(a+n)v\}$ for $n\in\frac12\Z.$
\end{lemm}
{\it Proof.}\ \ For any $a\in\C$, denote
$V(a)=\oplus_{n\in\frac12\Z}V^{a+n}$. From relations
(\ref{1eq1.1})--(\ref{1eq1.5}), one can easily see that $V(a)$ is an
$\LL[s]$-submodule of $V$ such that
$V=\oplus_{a\in\C/\frac12\Z}V(a)$ is a direct sum of different
$V(a)$. Hence $V=V(a)$ for some $a\in\C$.\hfill$\Box$\vskip2mm

\begin{lemm}\label{lemmsu1}
Suppose $V=\oplus_{a\in\C/\frac12\Z}V(a)$ is an irreducible
Harish-Chandra $\LL[s]$-module without highest and lowest weights.
Then for any $i\in\Z^*$, $k\in\frac12\Z$,
\begin{eqnarray}\label{s===0}
L_i|_{V_k}\oplus L_{i+1}|_{V_k}\oplus Y_{i+s}|_{V_k}\oplus
Y_{i+1+s}|_{V_k}: \ \ V_k\ \to\ V_{k+i}\oplus V_{k+i+1}\oplus
V_{k+i+s}\oplus V_{k+i+1+s}
\end{eqnarray}
is injective. In particular, by taking $i=-k$ $($if $k\in\Z)$ or
$i=\frac12-k$ $($if $k\in\frac12+\Z)$, we obtain that $\dim\,V_k$ is
uniformly bounded.
\end{lemm}{\it Proof.}\ \
Suppose there exists some $v_0\in V_k$ such that
\begin{equation}\label{LLL-1111}L_iv_0=L_{i+1}v_0=Y_{i+s}v_0=Y_{i+1+s}v_0=0.\end{equation}
Without loss of generality, we can suppose $i>0$. Note that when
$\ell\gg0$, we have
$$\ell=ai+b(i+1),\,\ \ \ell+s=a'i+b'(i+s),\,\ \ \ell=a''(i+s)+b''(i+1+s),$$
for some $a,b,a',b',a'',b''\in\N$, from this and (\ref{1eq1.1}),
(\ref{1eq1.2}) and (\ref{1eq1.4}), one can easily deduce that
$L_\ell,Y_{s+\ell},M_\ell$ can be generated by
$L_i,L_{i+1},Y_{i+s},Y_{i+1+s}$. Therefore there exists some $N>0$
such that
$$L_\ell v_0=Y_{s+\ell}v_0=M_\ell v_0=0\mbox{\ \ for all \ }\ell\ge
N.$$ The rest of the proof is exactly similar to that of
\cite[Proposition~2.1]{S5}.\QED\vskip2mm

Now Theorem \ref{maintheo}(i) follows from Lemma \ref{s===0}.

\vskip16pt

\lll{\bf3. \ Proofs of Theorem \ref{maintheo}(ii) and (iii)}
\setcounter{section}{3}\setcounter{theo}{0} \setcounter{equation}{0}

From now on, we shall suppose $V$ is a module of the intermediate
series over $\LL[s]$. As in \cite{S3,S4}, one sees that $c$ acts
trivially on $V$, so we can omit $c$ in (\ref{1eq1.1}).

First we prove Theorem \ref{maintheo}(iii), i.e., for the original
Schr\"{o}dinger-Virasoro Lie algebra $\LL[\frac12]$.

We shall first suppose that  ${\rm dim\,}V_k=1$ for all
$k\in\frac12\Z$ and both $V'=\oplus_{k\in\Z}V_k$ and
$V''=\oplus_{k\in\frac12+\Z}V_k$ are ${\cal V}ir$-modules of the
form $A_{a,b},\,a,b\in {\mathbb C}$ (later on, we shall determine
all possible deformations). Therefore, we can choose a basis
$\{x_k\,|\,k\in\frac{1}{2}\Z\}$ of $V$ such that
\begin{eqnarray}
Y_px_k\!\!\!&=&\!\!\!f_{p,k}x_{k+p},\label{eqactionfir2}\\
M_nx_k\!\!\!\!&=&\!\!\!g_{n,k}x_{k+n},\label{eqactionfir3}\\
L_nx_k\!\!&=&\!\!\!\left\{\begin{array}{ll}
\!\!(a+k+bn)x_{k+n}&{\rm if}\ \,k\in\Z,\vs{6pt}\\
\!\!(a+k+b'n)x_{k+n}&{\rm if}\ \,k\in\frac{1}{2}+\Z,
\end{array}\right.\label{eqactionfir11}
\end{eqnarray}
for some $a,\,b,\,b',\,f_{p,k},\,g_{n,k}\in{\C}$, where
$k\in\frac12\Z,\,p\in\frac{1}{2}+\Z,\,n\in\Z$.

\begin{lemm}\adddot\label{lemm2.2} $b'=b\pm\frac{1}{2}$ or $(b,b')
\in\big\{(0,\frac{3}{2}),\,(\frac{3}{2},0),\,(1,-\frac{1}{2}),\,(-\frac{1}{2},1)\big\}$.
\end{lemm}
{\it Proof.}\ \ If $f_{p,k}=0$ for all
$p\in\frac12+\Z,\,k\in\frac12\Z$, then from (\ref{1eq1.4}), we also
obtain $g_{n,k}=0$ for all $n\in\Z,\,k\in\frac12\Z$, and so $\SS$
acts trivially on $V$.  Thus there exists some
$p_0\in\frac{1}{2}+\Z$ and $k_0\in\frac{1}{2}\Z$ such that
$f_{p_0,k_0}\ne0$. Replacing $a$ by $a+\frac12$ if necessary (which
exchanges $V'$ and $V''$), we can always suppose $f_{p_0,k_0}\ne0$
for some $p_0\in\frac12+\Z,\,k_0\in\Z.$ Then from (\ref{1eq1.2}), we
see that for every $p\in \frac12+\Z$,
\begin{equation}\label{con-1}
f_{p,k}\ne0\mbox{ \ \ for infinitely many \ }k\in\Z.
\end{equation}
For any $m,\,n\in\Z,\,p\in\frac{1}{2}+\Z,\,k\in\Z$, applying
\begin{eqnarray}\label{eqactlly1}
&&(p-\frac{m+n}{2})[L_m,[L_n,Y_p]]=(p-\frac{n}{2})(n+p-\frac{m}{2})[L_{m+n},Y_p]
\end{eqnarray}
to $x_k$, using  (\ref{eqactionfir2}), (\ref{eqactionfir11}) and
comparing the coefficients of $x_{k+p+m+n}$, we obtain
\begin{eqnarray}\label{keyeqthree1}
\!\!\!\!\!\!\!\!\!\!\!\!\!\!\!\!\!\!\!\!\!\!\!\!
&&(p-\frac{m+n}{2})\Big((a+k+p+b'n)(a+k+p+n+b'm)f_{p,k}\nonumber\\[-2pt]
\!\!\!\!\!\!\!\!\!\!\!\!\!\!\!\!\!\!\!\!\!\!\!\!
&&\ \ \ \ \ \ \ \ \ \ \ -(a+k+bn)(a+k+p+n+b'm)f_{p,k+n}\nonumber\\[-2pt]
\!\!\!\!\!\!\!\!\!\!\!\!\!\!\!\!\!\!\!\!\!\!\!\! &&\ \ \ \ \ \ \ \ \
\ \ -(a+k+bm)(a+k+p+m+b'n)f_{p,k+m}
\nonumber\\[-2pt]\!\!\!\!\!\!\!\!\!\!\!\!\!\!\!\!\!\!\!\!\!\!\!\!&&\ \ \ \ \ \ \ \ \ \ \
+(a+k+bm)(a+k+m+bn)f_{p,k+m+n}\Big)\nonumber\\[-2pt]
\!\!\!\!\!\!\!\!\!\!\!\!\!\!\!\!\!\!\!\!\!\!\!\!
&&=(p-\frac{n}{2})(n+p-\frac{m}{2})\Big( (a+k+p+b'm+b'n)f_{p,k}
-(a+k+bm+bn)f_{p,k+m+n}\Big).
\end{eqnarray}
In the above equation, replacing $m,\,n,\,k$ by (i) $m,\,m,\,k-m$,
(ii) $-m,\,-m,\,k+m$ and (iii) $m,\,-m,\,k$, respectively, we obtain
the following three equations:
\begin{eqnarray}
\label{(I)}
&&\Big((p-\frac{m}{2})(p+\frac{m}{2})\big(a+k+(2b-1)m\big)\nonumber\\
&&\ \ \ \ \ \ \ \ \ \ \ +(p-m)\big(a+k+(b-1)m\big)(a+k+bm)\Big)f_{p,k+m}\nonumber\\
&&\ \ \ \ \ \ \ \ \ \ \
-2(p-m)\big(a+k+(b-1)m\big)(a+k+p+b'm)f_{p,k}\nonumber\\
&&\ \ \ \ \ \ \ \ \ \ \ +\Big((p-m)\big(a+k+p+(b'-1)m\big)(a+k+p+b'm)\nonumber\\
&&\ \ \ \ \ \ \ \ \ \ \
-(p-\frac{m}{2})(p+\frac{m}{2})\big(a+k+p+(2b'-1)m\big)\Big)f_{p,k-m}=0,\\
\label{(II)}
&&\Big((p+m)\big(a+k+p-(b'-1)m\big)(a+k+p-b'm)\nonumber\\
&&\ \ \ \ \ \ \ \ \ \ \ -(p+\frac{m}{2})(p-\frac{m}{2})\big(a+k+p-(2b'-1)m\big)\Big)f_{p,k+m}\nonumber\\
&&\ \ \ \ \ \ \ \ \ \ \ -2(p+m)\big(a+k-(b-1)m\big)(a+k+p-b'm)f_{p,k}\nonumber\\
&&\ \ \ \ \ \ \ \ \ \ \ +\Big((p+\frac{m}{2})(p-\frac{m}{2})\big(a+k-(2b-1)m\big)\nonumber\\
&&\ \ \ \ \ \ \ \ \ \ \ +(p+m)\big(a+k-(b-1)m\big)(a+k-b
m)\Big)f_{p,k-m}=0,
\end{eqnarray}
\begin{eqnarray}
\label{(III)}
&&(a+k+bm)\big(a+k+p-(b'-1)m\big)f_{p,k+m}\nonumber\\
&&\ \ \ \ \ \ \ \ \ \ \
-\Big((a+k+p-b'm)\big(a+k+p+(b'-1)m\big)\nonumber\\
&&\ \ \ \ \ \ \ \ \ \ \ \,+(p+\frac{m}{2})(\frac{3m}{2}-p)
+(a+k+bm)\big(a+k-(b-1)m\big)\Big)f_{p,k}\nonumber\\
&&\ \ \ \ \ \ \ \ \ \ \ +(a+k-b
m)\big(a+k+p+(b'-1)m\big)f_{p,k-m}=0.
\end{eqnarray}
Regard (\ref{(I)})--(\ref{(III)}) as a system of linear equations on
the unknown variables $f_{p,k+m},\,f_{p,k}$ and $f_{p,k-m}$. Denote
by $\D(p,k,m)$ the coefficient determinant. By (\ref{con-1}), we see
that for any pairs $(p,m)\in (\frac12+\Z)\times\Z$, there exist
infinitely many integers $k$ such that $\D(p,k,m)=0$. Since
$\D(p,k,m)$ is a polynomial on $p,\,k,\,m$, we obtain that
\begin{equation} \label{Delta=0}
\D(p,k,m)=0 \mbox{ \ for all \ }p\in\frac{1}{2}+\Z,\ k,\,m\in\Z.
\end{equation}
It is a little lengthy but straightforward to compute (one can
simply use {\it Mathematica} to solve a system of linear equations
without problem)
\begin{eqnarray*}
\D(p,k,m)\!\!\!&=&\!\!\!\frac{m^6}{64}\D_0(\D_1(a+k)p+\D_2m^2+\D_3p^2),
\end{eqnarray*}\vs{-6pt}
where\vs{2pt}
\begin{eqnarray*}
\D_0\!\!\!&=&\!\!\!(2b-2b'-1)(2b-2b'+1),\\
\D_1\!\!\!&=&\!\!\!18(2b+2b'-3)(2b+2b'-1),\\
\D_2\!\!\!&=&\!\!\!
4(b+b'-1)(-3b-4b^2+4b^3-9b'+4bb'+4b^2b'+12{b'}^2-4b{b'}^2-4{b'}^3),
\\
\D_3\!\!\!&=&\!\!\!27-156b+152b^2-16b^3-16b^4-180b'+328bb'-80b^2b'\\
&&-32b^3b'+240{b'}^2+240{b'}^2-176b{b'}^2-112{b'}^3+32b{b'}^3+16{b'}^4.
\end{eqnarray*}
It is easy to see that (\ref{Delta=0}) holds if and only if
\begin{eqnarray*} \label{Delta=022}
&&\D_0=0\ \ \mbox{or}\ \ \D_1=\D_2=\D_3=0,
\end{eqnarray*}
if and only if
\begin{eqnarray}\label{eqrelbb'}
&&b'=b\pm \frac{1}{2}\ \ \mbox{or}\ \
(b,\,b')\in\big\{(0,\,\frac{3}{2}),\ \,(\frac{3}{2},\,0),\
\,(1,\,-\frac{1}{2})\ \,\mbox{or}\ \,(-\frac{1}{2},\,1)\big\}.
\end{eqnarray}
(Note that from $\D_1=0$, one obtains that $b'=\frac32-b$ or
$\frac12-b$, then from $\D_2=\D_3=0$ one can solve $b$.) Thus this
lemma follows.\QED\vskip2mm

By replacing $a$ by $a+\frac12$ (then (\ref{con-1}) does not
necessarily  hold) if necessary (which exchanges $V'$ and $V''$), we
can always suppose ${\rm Re}(b')\ge{\rm Re}(b)$ (where ${\rm Re}(b)$
is the real part of the complex number $b$). Thus we only need to
consider the cases
\begin{eqnarray}\label{b-b1}
b'=b+\frac{1}{2}\mbox{ \ \ or \ }(b,b')
\in\big\{(0,\frac{3}{2}),\,(-\frac{1}{2},1)\big\}.
\end{eqnarray}
We shall now consider all the cases given in (\ref{b-b1}) one by
one.
\begin{case}\adddot\label{case1}
$b'=b+\frac{1}{2}$.
\end{case}

We need to determine $f_{p,k},\,g_{n,k}$ defined in
(\ref{eqactionfir2}) and (\ref{eqactionfir3}), where $n\in\Z$,
$p\in\frac{1}{2}+\Z$ and $k\in\frac{1}{2}\Z$. From the relation
$[Y_m,Y_p]x_k=(p-m)M_{m+p}x_k$ for
$m,p\in\frac{1}{2}+\Z,\,k\in\frac{1}{2}\Z$, we only need to
determine $f_{p,k}$.
\begin{lemm}\adddot\label{Lemma-case1}
For any $m,\,k\in\Z,\,p\in\frac{1}{2}+\Z$,
\begin{eqnarray}
&&(a+k+2bp)f_{p,k+m}=(a+k+m+2bp)f_{p,k}.\label{eqfpk+mpkI11}
\end{eqnarray}
\end{lemm}
{\it Proof.}\ \ Using (\ref{(I)}) and (\ref{(II)}) to cancel
$f_{p,k-m}$, one can obtain (or simply using {\it Mathematica}) that
(\ref{eqfpk+mpkI11}) holds under the following condition
\begin{eqnarray}\label{eqconddelkmp1}
\Delta_1(k,m)\!\!\!&=&\!\!\!\Big((a+k)(1-4b)+(1+6b)(b-1)p\Big)m^2\nonumber\\
&&+6(a+k)^2p+8(a+k)(1-b)p^2+4(1-b)p^3\neq0.
\end{eqnarray}
Next we want to prove that (\ref{eqfpk+mpkI11}) holds  under the
condition
\begin{eqnarray}\label{eqconddelkmp1+1}
\Delta_2(k)\!\!\!&=&\!\!\!(a+k)(1-4b)+(1+6b)(b-1)p \neq0.
\end{eqnarray}
Suppose (\ref{eqconddelkmp1+1}) holds. Then $\D_1(k+m',m-m')$ is a
polynomial on $m'$ of degree 3 (if $b\ne\frac14$) or of degree 2 (if
$b=\frac14$, in this case the coefficient of $m'^2$ in
$\D_1(k+m',m-m')$ is $(1+6b)(b-1)p+6p=(6-\frac{5}{8})p\ne0$ since
$p\in\frac12+\Z$), and $\D_1(k,m')$ is a polynomial on $m'$ of
degree $2$. Thus we can find some $m'$ such that
$$(a+k+m'+2bp)\D_1(k+m',m-m')\D_1(k,m')\ne0,$$ in particular,
condition (\ref{eqconddelkmp1}) holds for the two triples
$(k+m',m-m',p),\,(k,m',p)$.
 Thus (\ref{eqfpk+mpkI11}) holds for these two triples, and one has
\begin{eqnarray}\label{eqfpk+mpkI211}
(a+k+2bp)f_{p,k+m}\!\!\!&=&\!\!\!\frac{a+k+2bp}{a+k+m'+2bp}(a+(k+m')+2bp)f_{p,(k+m')+(m-m')}\nonumber\\
\!\!\!&=&\!\!\!\frac{a+k+2bp}{a+k+m'+2bp}(a+(k+m')+(m-m')+2bp)f_{p,k+m'}\nonumber\\
\!\!\!&=&\!\!\!\frac{a+k+m+2bp}{a+k+m'+2bp}(a+k+2bp)f_{p,k+m'}\nonumber\\
\!\!\!&=&\!\!\!(a+k+m+2bp)f_{p,k}.
\end{eqnarray}
Thus (\ref{eqfpk+mpkI11}) holds under condition
(\ref{eqconddelkmp1+1}). Now for any $p,k,m$, choose $m'$ such that
\begin{equation}\label{3223}
(a+k+m'+2bp)\D_2(k+m')\ne0.\end{equation}
 This follows by noting
that if $b\ne\frac14$ then $\D_2(k+m')$ is a polynomial on $m'$ of
degree 1, and if $b=\frac14$ then
$\D_2(k+m')=(1+6b)(b-1)p=-\frac{15}{8}p\ne0$ (since
$p\in\frac12+\Z$). Now (\ref{3223}) implies that we have all
equalities of (\ref{eqfpk+mpkI211}) except the last equality. But
the last equality also follows by writing $f_{p,k}$ as
$f_{p,(k+m')+(-m')}$ and using (\ref{eqfpk+mpkI11}) with the pair
$(k,m)$ being replaced by $(k+m',-m')$ (condition
(\ref{eqconddelkmp1+1}) holds for this pair). \QED

For any $m\in\Z,\,p\in\frac{1}{2}+\Z,\,k\in\frac{1}{2}\Z$, applying
$[L_m,Y_p]=(p-\frac{m}{2})Y_{p+m}$ to $x_k$ and comparing the
coefficients of $x_{k+m+p}$, one has
\begin{eqnarray}\label{eqfnkfnm1}
\left\{\begin{array}{cc}\!\!\!\!\!\!\!\!\!\!\!\!
(a+k+bm)f_{p,m+k}-(a+k+p+b'm)f_{p,k}
=\frac{m-2p}{2}f_{m+p,k}\,\ \!\ \ {\rm if}\,k\in\Z,\vs{6pt}\\
\!\!(a+k+b'm)f_{p,m+k}-(a+k+p+bm)f_{p,k}=\frac{m-2p}{2}f_{m+p,k}\
\,\ \ {\rm if}\,k\in\frac{1}{2}+\Z.
\end{array}\right.
\end{eqnarray}
Taking $m=p-n$ in (\ref{eqfnkfnm1}) and using $b'=b+\frac{1}{2}$,
one has
\begin{eqnarray}\label{eqfnk??m1??}
(a+k+bp-b n)f_{n,p-n+k}-(a+k+n+bp-b n+\frac{p-n}{2})f_{n,k}
=\frac{p-3n}{2}f_{p,k}.
\end{eqnarray}
Using (\ref{eqfpk+mpkI11}), one has
\begin{eqnarray}\label{eqaffnk1}
&&(a+k+2b n)f_{n,k+p-n}=(a+k+p-n+2b n)f_{n,k}.
\end{eqnarray}
Then using (\ref{eqfnk??m1??}) together with (\ref{eqfpk+mpkI11})
and (\ref{eqaffnk1}), we obtain
\begin{eqnarray}\label{eqfk+mfnk1}
&&(p-3n)(a+k+2b p)\big((a+k+2b n)f_{p,k+m}-(a+k+m+2b
p)f_{n,k}\big)=0.
\end{eqnarray}
Letting $m=j-k$ in (\ref{eqfk+mfnk1}), one has
\begin{eqnarray}\label{eqfk+mfnk2}
&&(p-3n)(a+k+2b p)\big((a+k+2b n)f_{p,j}-(a+j+2b p)f_{n,k}\big)=0.
\end{eqnarray}
\begin{lemm}\adddot\label{eqfclaim11}
For any $k,\,j\in\Z,\,n,\,p\in\frac{1}{2}+\Z$, one has
\begin{eqnarray}\label{eqcaim11}
&&(a+k+2b n)f_{p,j}=(a+j+2b p)f_{n,k}.
\end{eqnarray}
Furthermore, for any $j\in\Z,\,p\in\frac{1}{2}+\Z$, $f_{p,j}$ can be
written as
\begin{eqnarray}\label{eqcaim12}
&&f_{p,j}=(a+j+2b p)f_0,
\end{eqnarray}
for some $f_0\in\C$.
\end{lemm}
{\it Proof.} Equation (\ref{eqcaim12}) follows from (\ref{eqcaim11})
by fixing $n_0,k_0$ such that $a+k_0+bn_0\ne0$ and letting
$f_0=\frac{f_{n_0,k_0}}{a+k_0+bn_0}$. It remains to prove
(\ref{eqcaim11}). If $(p-3n)(a+k+2b p)\neq0$, then (\ref{eqcaim11})
follows from (\ref{eqfk+mfnk2}). Next suppose $(p-3n)(a+k+2b p)=0$.
Choose $n_0\in\frac{1}{2}+\Z$ and $k_0\in\Z$ such that
$$(p-3n_0)(a+k_0+2b p)\neq0\ \mbox{ \ and \ }\ (n-3n_0)(a+k_0+2b n)\neq0.$$ Then by
(\ref{eqfk+mfnk2}), $(a+k_0+2b n_0)f_{p,j}=(a+j+2b p)f_{n_0,k_0}$
and also $(a+k+2b n)f_{n_0,k_0}=(a+k_0+2b n_0)f_{n,k}$. Thus
(\ref{eqcaim11}) also holds. This proves the lemma.\QED\vskip2mm

As for the case $k\in\frac{1}{2}+\Z$, we have the following lemma.
\begin{lemm}\adddot\label{Lemma2-case1}
For any $m,k,p\in\frac{1}{2}+\Z$,
\begin{eqnarray}
&&f_{p,k}=d_{0}.\label{eqfpkfpkin12Z}
\end{eqnarray}
\end{lemm}
{\it Proof.} First suppose $b\ne0$. We want to prove
\begin{equation}\label{Su0}
f_{\frac12,k}=d_0\mbox{ \ \ for \ }k\in\frac12+\Z.
\end{equation}
Applying $[L_1,Y_{\frac12}]=0$ to $x_k$ for $k\in\frac12+\Z$ and
using $b'=b+\frac12$, we obtain \begin{equation}\label{Su1}
(a+k+b+\frac12)f_{\frac12,k}-f_{\frac12,k+1}(a+k+b+\frac12)=0.
\end{equation}
If $a+b\notin\Z$, then $a+k+b+\frac12\ne0$ for all $k\in\frac12+\Z$,
thus (\ref{Su1}) implies $f_{\frac12,k}$ is a constant, denoted by
$d_0$. Thus (\ref{Su0}) holds in this case. Assume that $a+b\in\Z$.
If necessary, by shifting the index $k$ of $x_k$ by an integer
(which does not change $V',\,V''$), we can suppose $a+b=0$. Then
(\ref{Su1}) shows
\begin{equation}\label{Su2}
f_{\frac12,k}=\left\{\begin{array}{ll} f_{\frac12,\frac12}&\mbox{if
\ }k>0,\\[6pt]
f_{\frac12,-\frac12}&\mbox{if \
}k<0.\end{array}\right.\end{equation} Applying
$[L_{-1},Y_{\frac12}]=Y_{-\frac12}$ to $x_k$, we obtain
$f_{-\frac12,k}=(a+k+\frac12-b)f_{\frac12,k}-f_{\frac12,k-1}(a+k-(b+\frac12))$,
which together with (\ref{Su2}) gives (using $a=-b$)
\begin{equation}\label{Su3}
f_{-\frac12,k} =\left\{\begin{array}{clc}
f_{\frac12,\frac12}&\mbox{if \
}k>\frac12,\\[6pt]
(1-2b) f_{\frac12,\frac12}+2bf_{\frac12,-\frac12}&\mbox{if \
}k=\frac12,\\[6pt]
 f_{\frac12,-\frac12}&\mbox{if \
}k\le-\frac12.
\end{array}
\right.\end{equation} Applying $[L_1,Y_{-\frac12}]=-Y_{\frac12}$ to
$x_{\frac12}$, using (\ref{Su2})  and (\ref{Su3}),  we obtain
$$(1-2b) f_{\frac12,\frac12}+2bf_{\frac12,-\frac12}=
-f_{\frac12,\frac12}=(a+\frac12-\frac12+b)f_{-\frac12,\frac12}-f_{-\frac12,\frac32}(a+\frac12+b+\frac12),$$
which implies $f_{\frac12,\frac12}=f_{\frac12,-\frac12}$ since
$a=-b\ne0$. This together with (\ref{Su2}) gives (\ref{Su0}).

Now for any $\frac23\ne p\in\frac12+\Z$, by replacing $(n,p)$ by
$(p-\frac12,\frac12)$ in (\ref{1eq1.2}) and applying it to $x_k$,
one can easily obtain $f_{p,k}=d_0$. If $p=\frac23$, choosing some
$n\in\Z,\,p\in\frac12+\Z$ such that $n+p=\frac23$ and
$p-\frac{n}{2}\ne0$, we can again obtain $f_{\frac32,k}=d_0$. This
proves (\ref{eqfpkfpkin12Z}) for the case $b\ne0$.

Now suppose $b=0$. Similar to the proof of (\ref{eqfpk+mpkI11}), for
any $m,n\in\Z,\,k,p\in\frac{1}{2}+\Z$, one has
\begin{eqnarray}
&&f_{p,k+m}=f_{p,k},\label{eqfpkm123}
\end{eqnarray}
under the condition (noting that $b=0$)
\begin{eqnarray}\label{N1}
&&\nabla(k,m)=\nabla_2(k) m^2+\nabla_0(k)\neq0,
\end{eqnarray}
where
\begin{eqnarray*}
\nabla_2(k)\!\!\!&=&\!\!\!(a+k)^2+3(a+k)p+2p^2,\\
\nabla_0(k)\!\!\!&=&\!\!\!-6(a+k)^3p-10(a+k)^2p^2-6(a+k)p^3-2p^4.
\end{eqnarray*}
From (\ref{N1}), one can easily find some $m'$ such that (since the
left-hand side of (\ref{nakm1m}) is a polynomial on $m'$ of degree
$8$)
\begin{eqnarray}\label{nakm1m}
\nabla(k+m',m-m')\nabla(k+m',-m')\ne0.
\end{eqnarray}
Then (\ref{eqfpkm123}) holds for the pairs
$(k+m',m-m'),\,(k+m',-m')$, which forces
\begin{eqnarray*}
&&f_{p,k+m}=f_{p,(k+m')+(m-m')}=f_{p,(k+m')}=f_{p,(k+m')+(-m')}=f_{p,k}.
\end{eqnarray*}
Thus the lemma follows.\QED\vskip2mm

For any $m,\,n\in\frac{1}{2}+\Z,\,k\in\Z$, applying
$[Y_m,Y_n]=(n-m)M_{m+n}$ to $x_k$ and comparing the coefficients of
$x_{k+m+n}$, one has
\begin{eqnarray}\label{eqggm1}
&&(n-m)(g_{m+n,k}-2bd_0f_0)=0.
\end{eqnarray}
Taking $n=p-m$ with $p-2m\ne0$ in (\ref{eqggm1}), one immediately
obtains
\begin{eqnarray}\label{eqfggm122}
&&g_{p,k}=2bd_0f_0\ \ \ \mbox{for all}\ k,\,p\in\Z.
\end{eqnarray}
Now applying $[Y_m\,,Y_n]=(n-m)M_{m+n}$ to $x_k$ for
$k\in\frac{1}{2}+\Z$, and comparing the coefficients of $x_{k+m+n}$,
one has
\begin{eqnarray}\label{eqfnkfgm1}
&&(n-m)\big(g_{m+n,k}-(1-2b)d_0f_0\big)=0.
\end{eqnarray}
Thus as in (\ref{eqfggm122}), we obtain
\begin{eqnarray}\label{eqfggm122+}
&&g_{p,k}=(1-2b)d_0f_0\ \ \ \mbox{for all}\
p\in\Z,\,k\in\frac{1}{2}+\Z.
\end{eqnarray}
By now, we have obtained the following relations
\begin{eqnarray}\label{eqfg36}
\!\!\!\!\!\!\!\!\!\!\!&&f_{p,k}=\left\{\begin{array}{cc}\!\!\!\ \ \
\ \ \ \ \ \ d_0\ \,\ \,\
\ \ \ \ \,\ \ {\rm if}\,\,k\in\frac{1}{2}+\Z,\vs{6pt}\\
\!\!\!\!\!\!\!\!\!\!\!\!(a+k+2b p)f_0\,\ {\rm if}\,\,k\in\Z,
\end{array}\right.\
g_{n,k}=\left\{\begin{array}{cc}(1-2b)d_0f_0\,\,\ \ {\rm if}\,\,k\in\frac{1}{2}+\Z,\vs{6pt}\\
\!\!\!\!\!\!\,2bd_0f_0\,\ \ \ \ \ \ \ \ {\rm if}\,\,k\in\Z.
\end{array}\right.
\end{eqnarray}

For any $p\in\frac{1}{2}+\Z,\,n,k\in\Z$, applying $[Y_p,M_n]=0$ to
$x_k$ and comparing the coefficients of $x_{k+n+p}$, one has
\begin{eqnarray}\label{eqbd_0f01}
&&2(a+n+k+2b p)bd_0f_0^2=(a+k+2b p)(1-2b)d_0f_0^2.
\end{eqnarray}
Comparing the coefficients of $n^1$ and $n^0$, we obtain
$bd_0f_0^2=(1-2b)d_0f_0^2=0$, which implies
\begin{eqnarray}\label{eqbd_0f023}
&&d_0f_0=0.
\end{eqnarray}
Since $V$ is indecomposable, we obtain
\begin{eqnarray}\label{eqbd_0f02368}
&&d_0=0,\ \,f_0\ne0\ \,{\rm{\ or \ }}\ \,d_0\ne0,\ \,f_0=0,
\end{eqnarray}
and the second equation  of (\ref{eqfg36}) can be rewritten as
\begin{eqnarray}\label{eqbd_0f0236}
&&g_{n,k}=0\ \ \ \mbox{for all}\ \,n\in\Z,\,k\in\frac{1}{2}\Z.
\end{eqnarray}
For the first case of  (\ref{eqbd_0f02368}), by rescaling $x_k$ for
$k\in\Z$ (and keeping $x_k,\,k\in\frac12+\Z$ unchanged), we can
suppose
\begin{eqnarray}\label{eqbd_0f0236822}
&&f_0=1.
\end{eqnarray}
Similarly, for the last case of  (\ref{eqbd_0f02368}), we can
suppose
\begin{eqnarray} \label{eqbd_0fd0=1}
&&d_0=1.
\end{eqnarray}
Then we get the two types of irreducible modules of intermediate
series over $\LL$, denoted by $A_{a,b}$, $B_{a,b}$ with the  basis
$\{x_k\,|\, k\in\frac{1}{2}\Z\}$ and the actions given by
(\ref{caseIIact1})--(\ref{caseIIact3}) and
(\ref{caseIIIac1})--(\ref{caseIIIac3}) respectively.

Now we consider all possible deformations of the representation
$A_{a,b}$. In order for $A_{a,b}$ to have a nontrivial deformation,
it is necessary that $A_{a,b}$ is decomposable. Thus we only need to
consider the following 4 subcases.

\noindent{\bf Subcase 1.1}\quad$a\in\Z,\,\,b=1$.

In this case, by shifting the index $k$, we can suppose $a=0,\,b=1$.
Then $x_k,\, 0\ne k\in\frac12\Z$ span the only proper submodule.
Thus the possible deformations are the actions of $M_n,Y_p,L_n$ on
$x_0$.

In order to get nontrivial deformation, we firstly remark that the
actions $L_n$ on $x_0$ must have deformation, otherwise from our
computations above, the actions of $M_n,Y_p$ would not have any
deformation. Furthermore, the deformation of the action of $L_n$ on
$x_0$ can be seen in (\ref{condition++}). Thus according to
(\ref{caseIIact1})--(\ref{caseIIact3}), we have
\begin{eqnarray}
M_nx_k\!\!\!\!&=&\!\!\!0\ \ \ \,\mbox{for any}\
k\in\frac{1}{2}\Z^*,\label{2caseIIact1}\vs{6pt}\\
Y_px_k\!\!\!&=&\!\!\!\left\{\begin{array}{cl}0&{\rm if}\ \,k\in\frac{1}{2}+\Z,\vs{6pt}\\
(k+2p)x_{k+p}&{\rm if}\ \,k\in\Z^*,
\end{array}\right.\label{2caseIIact2}\\
L_nx_k\!\!&=&\!\!\!\left\{\begin{array}{cc}\!(k+\frac{3}{2}n)x_{k+n}\
\,
\,{\rm if}\ \,k\in\frac{1}{2}+\Z,\vs{6pt}\\
\!\!\!\!\!\!\!\!\!\!\!\!\!\!\ \ \ (k+n)x_{k+n}\! \ \, \,\ \ {\rm
if}\ \,k\in\Z^*, \vs{6pt}\\ \!\!\!\!\!\!\!\!\!\!n(n+\a)x_{n}\ \, \,\
\ {\rm if}\ \,k=0,
\end{array}\right.\label{2caseIIact3}
\end{eqnarray}
for some $\a\in\C$,  where $p\in\frac{1}{2}+\Z,\,n\in\Z$. And we
suppose
\begin{equation}\label{SSS}
M_n x_0=f_n x_n,\ \ Y_p x_0=g_p x_p ,\ \ n\in\Z,\
p\in\frac12+\Z,\end{equation} for some $f_n,\,g_p,\,h_n\in\C$. For
any $n\in\Z^*,\,p\in\frac{1}{2}+\Z$, applying
$[L_{n},Y_p]=(p-\frac{n}{2})Y_{n+p}$ to $x_{0}$ and comparing the
coefficients of $x_{n+p}$, we obtain
\begin{eqnarray}\label{eqfg8fs3666}
&&(p+\frac{3n}{2})g_p-n(n+\alpha)(n+2p)=(p-\frac{n}{2})g_{n+p}.
\end{eqnarray}
Letting $n=2p$ in (\ref{eqfg8fs3666}), we obtain
$g_p=2p(2p+\alpha),$ and also $g_{n+p}=2(n+p)(2(n+p)+\alpha).$ Using
them in (\ref{eqfg8fs3666}), we obtain a contradiction, i.e., in
this case $A_{a,b}$ has no deformation.

\noindent{\bf Subcase 1.2}\quad$a\in\Z,\,\,b=0$.

Similarly, we can suppose $a=b=0$ by shifting the index $k$. Then
$\F x_0$ is the only proper submodule (thus the deformation of the
action of $L_n$ on $x_{-n}$ can be seen in (\ref{condition+++})).
Hence in order to obtain the identities listed in
(\ref{2caseIIsub2act1})--(\ref{2caseIIsub2act3}), we only need to
compute the complex numbers $g_p$ and $f_n$ defined by
\begin{eqnarray}
&&Y_px_{-p}=g_px_{0},\ \ \,M_nx_{-n}=f_nx_0 \ \ \mbox{ for all}\ \
n\in{\ma Z},\,\,p\in\frac{1}{2}+\Z.
\end{eqnarray}
For any $n\in\Z^*,\,p\in\frac{1}{2}+\Z$, applying
$[L_{n},Y_p]=(p-\frac{n}{2})Y_{n+p}$ to $x_{-n-p}$ and comparing the
coefficients of $x_{0}$, we obtain $g_p=0$ for any
$p\in\frac{1}{2}+\Z$. Then
\begin{eqnarray}\label{eqfgp=0CaseII}
&&Y_px_{-p}=0\ \ \mbox{for any}\ \,p\in\frac{1}{2}+\Z.
\end{eqnarray}
By (\ref{eqfgp=0CaseII}) and using
$[Y_p,Y_{n-p}]x_{-n}=(n-2p)M_{n}x_{-n}$, one has
\begin{eqnarray}\label{eqfmn-n=0CaseII}
&&M_nx_{-n}=0\ \ \mbox{for any}\ \,n\in\Z.
\end{eqnarray}
Then we get a deformation of $A_{a, b}$, denoted by $A_2(\a)$ with
the relations given by
(\ref{2caseIIsub2act1})--(\ref{2caseIIsub2act3}).

\noindent{\bf Subcase 1.3}\quad
$a\in\frac{1}{2}+\Z,\,\,b=\frac{1}{2}$.

Similarly, shifting the index $k$, one can assume
$a=-\frac{1}{2},\,b=\frac{1}{2}$. Then $x_k,\frac12\ne
k\in\frac12\Z$, span the only submodule. Using the same techniques,
we can obtain
\begin{eqnarray}\label{eqCaseII.3YM01}
&&M_nx_\frac{1}{2}=Y_px_\frac{1}{2}=0\ \ \mbox{for
all}\,\,n\in\Z,\,p\in\frac{1}{2}+\Z.
\end{eqnarray}
Then we can obtain a deformation, denoted by $A_1(\a)$ with the
relations given by (\ref{3caseIIact122})--(\ref{3caseIIact322}) (in
fact $A_1(\a)$ is isomorphic to the dual module of $A_2(\a')$ for
some $\a'\in\C$).

\noindent{\bf Subcase 1.4}
\quad$a\in\frac{1}{2}+\Z,\,\,b=-\frac{1}{2}$.

In this case, by shifting the index $k$, we can suppose
$a=-\frac12,\,b=-\frac12$. Then by
(\ref{caseIIact1})--(\ref{caseIIact3}), one has
\begin{eqnarray}
M_nx_k\!\!\!&=&\!\!\!0\ \ \ \,{\rm for\ any}\ k\in\frac{1}{2}\Z,\,n\neq\frac{1}{2}-n,\label{caseIIact122}\\
Y_px_k\!\!\!&=&\!\!\!\left\{\begin{array}{ll}
\ \ \ \ \ \ 0&{\rm if}\ \,k\in\frac{1}{2}+\Z,\,k\neq\frac{1}{2}-p,\vs{3pt}\\
(-\frac{1}{2}+k-p)x_{k+p}&{\rm if}\ \,k\in\Z,
\end{array}\right.\label{caseIIact222}\\
L_nx_k\!\!\!&=&\!\!\!\left\{\begin{array}{lll}(-\frac{1}{2}+k)x_{k+n}&{\rm
if}\ \,
k\in\frac{1}{2}+\Z,\,k\ne\frac{1}{2}-n,\vs{3pt}\\
(-\frac{1}{2}+k-\frac{1}{2}n)x_{k+n}&{\rm if}\ \,\,k\in\Z,\vs{3pt}\\
-n(n+\a)x_{\frac{1}{2}}&{\rm if}\ \,k=\frac{1}{2}-n,
\end{array}\right.\label{caseIIact322}
\end{eqnarray}
where $p\in\frac{1}{2}+\Z,\,n\in\Z,\,\a\in\C$.

In order to determine all possible deformations, we suppose
\begin{eqnarray}\label{eqgp00}
&&M_nx_{\frac{1}{2}-n}=f_nx_{\frac{1}{2}},\ \ \
Y_px_{\frac{1}{2}-p}=g_px_{\frac{1}{2}} \ \ \mbox{ for all}\
\,n\in\Z,\,p\in\frac{1}{2}+\Z.
\end{eqnarray}
For any $n\in\Z^*,\,p\in\frac{1}{2}+\Z$, applying
$[L_{n},Y_p]=(p-\frac{n}{2})Y_{n+p}$ to $x_{\frac{1}{2}-n-p}$ and
comparing the coefficients of $x_{\frac{1}{2}}$, we obtain
\begin{eqnarray}\label{eqsbc1.41}
&&(p+\frac{3n}{2})g_p+n(n+2p)(n+\a)=(p-\frac{n}{2})g_{n+p}.
\end{eqnarray}
Letting $n=2p$ in (\ref{eqsbc1.41}), we obtain $g_p=-2p(2p+\a)$. As
in  Subcase 1.1, we obtain a contradiction with (\ref{eqsbc1.41}).
So in this case $A_{a,b}$ has no deformation either. \vskip2mm

 Similarly, we can also give all possible deformations of
the representation $B_{a,b}$ denoted by $B_1(\a)$ and $B_2(\a)$,
whose module structures are listed in
(\ref{caseIIIac1a})--(\ref{caseIIIac3b}).
\begin{case}\adddot\label{case2}
 $b=0,\ \,b'=\frac{3}{2}$.
\end{case}

Regardless of their deformations, the module structures over $\LL$
corresponding to $(b,b')=(0,\frac32)$ can be determined as those
 corresponding to the case $(b,b')=(1,\frac{3}{2})$,
i.e., $A_{1,\frac{3}{2}}$ and $B_{1,\frac{3}{2}}$, whose module
structures have been listed in
(\ref{caseIIact1})--(\ref{caseIIIac3}) explicitly. Using the same
techniques, we obtain an indecomposable module of intermediate
series over $\LL$, denoted by $C_a$, given by
(\ref{eqaCa1})--(\ref{eqaCa3}) and its possible deformation, denoted
by $C(\a,\a')$ and given by (\ref{eqacIIS1+})--(\ref{eqacIIS3}).

\begin{case}\adddot\label{case2+}
$b=-\frac{1}{2},\ \,b'=1$.
\end{case}

Regardless of their deformations, the module structures over $\LL$
corresponding to $(b,b')=(-\frac{1}{2},1)$ can be determined as
those corresponding to the case $(b,b')=(-\frac{1}{2},0)$, i.e.,
$A_{-\frac{1}{2},0}$ and $B_{-\frac{1}{2},0}$, whose module
structures have been listed in
(\ref{caseIIact1})--(\ref{caseIIIac3}) explicitly. Using the same
techniques, we obtain an indecomposable module of intermediate
series over $\LL$, denoted by $D_a$, given by
(\ref{eqaDa1})--(\ref{eqaDa3}) and its  possible deformation,
denoted by $D(\b,\b')$ and given by
(\ref{eqLLSSI1})--(\ref{eqLLSSI3}).

This completes the proof of Theorem \ref{maintheo}(iii).

Finally we consider the case $\LL[0]$. In this case, $V_k=0$ for
$k\in\frac12+\Z$, and so $V$ is $\Z$-graded. As the proof above (now
all indices $k,p,n$ are in $\Z$), we still can suppose
(\ref{eqactionfir2}), (\ref{eqactionfir3}) and the first case of
(\ref{eqactionfir11}). Thus we again have (\ref{keyeqthree1}) with
$b'$ being $b$. But the proof of Lemma \ref{lemm2.2} shows that
$b'=b$ is impossible, i.e., $\SS$ must act trivially on $V$. This
proves Theorem \ref{maintheo}(ii).\QED

\end{document}